\documentclass[aos,MSNbibl,seceqn,nameyear,dvips]{arximspdf}

\doi{10.1214/11-AOS941}
\volume{40}
\issue{1}
\pubyear{2012}
\firstpage{30}
\lastpage{44}

\makeatletter

\newtheorem{theorem}{Theorem}[section]
\newtheorem{lemma}{Lemma}[section]

\newcommand{\boldT}{\mathbf{T}}
\newcommand{\boldX}{\mathbf{X}}

\makeatother

\begin{document}
\begin{frontmatter}

\title{Sequential monitoring with conditional randomization tests}
\runtitle{Sequential monitoring of conditional randomization tests}

\begin{aug}
\author[A]{\fnms{Victoria} \snm{Plamadeala}\ead[label=e1]{vplamadeala@ptilabs.com}}
\and
\author[B]{\fnms{William F.} \snm{Rosenberger}\corref{}\thanksref{T1}\ead[label=e3]{wrosenbe@gmu.edu}\ead[label=u1,url]{http://www.foo.com}}
\runauthor{V. Plamadeala and W. F. Rosenberger}
\affiliation{Precision Therapeutics and George Mason University}
\address[A]{Precision Therapeutics\\
2516 Jane Street\\
Pittsburgh, Pennsylvania 15203\\
USA\\
\printead{e1}}
\address[B]{Department of Statistics\\
George Mason University\\
4400 University Drive, MS 4A7\\
Fairfax, Virginia 22030\\
USA\\
\printead{e3}} %adresu isvedimo komanda gale!
\end{aug}

\thankstext{T1}{Supported by NSF Grant DMS-09-04253
under the 2009 American Reinvestment and Recovery Act.}

% HISTORY:
\received{\smonth{7} \syear{2011}}
\revised{\smonth{11} \syear{2011}}

% ABSTRACT
%
\begin{abstract}
Sequential monitoring in clinical trials is often employed to allow for
early stopping and other interim decisions, while maintaining the type
I error rate. However, sequential monitoring is typically described
only in the context of a population model. We describe a computational
method to implement sequential monitoring in a~randomization-based
context. In particular, we discuss a new technique for the computation
of approximate conditional tests following restricted randomization
procedures and then apply this technique to approximate the joint
distribution of sequentially computed conditional randomization tests.
We also describe the computation of a randomization-based analog of the
information fraction. We apply these techniques to a restricted
randomization procedure, Efron's [\textit{Biometrika} \textbf{58}
(1971) 403--417] biased coin design. These techniques require
derivation of certain conditional probabilities and conditional
covariances of the randomization procedure. We employ combinatoric
techniques to derive these for the biased coin design.
\end{abstract}

% KEYWORDS
%

\end{frontmatter}

%s1 #&#
\section{Introduction}

Sequential monitoring refers to analyzing data periodically during the
course of a clinical trial, with the purpose of detecting early
evidence in support of or against a hypothesis. A desirable feature of
such a monitoring plan would be flexible inspections of the data that
can occur at arbitrary time points. At the same time, sequentially
tested hypotheses must maintain the overall probability of type I error
at the prespecified level, since repeated testing is known to inflate
it. The Lan and DeMets (\citeyear{LanDeM83}) error spending approach
for sequential monitoring allows this. The approach makes use of a type
I error spending function, which depends on the amount of ``statistical
information'' available at the time of the interim inspection. In the
context of sequential monitoring, the statistical information is a
measure of how far a trial has progressed. Under a population model,
the amount of interim information---the information fraction---is
defined as the proportion of Fisher's information observed thus far in
the trial. The type I error spending function rations the amount of
type I error that may be spent at each look commensurate to the
information fraction. The critical value associated with the allowable
probability of type I error at a certain interim look is obtained and
compared to the observed value of the statistic. The decision whether
to continue, or stop, the trial is based on this comparison. Sequential
monitoring is typically discussed in the context of a population model.
However, it is not uncommon for the FDA to require a ``re-analysis'' of
data using a ``re-randomization test,'' or, as we call it here, a
\textit{randomization test}, defined below.

Let
$\boldT=T_1,\ldots,T_n$ be a randomization sequence, where $T_i=1$ if
subject $i$ is randomized to treatment 1; $T_i=0$ if subject $i$ is
randomized to treatment 2, $i=1,\ldots,n$. Let $N_1(j)=\sum_{i=1}^j
T_i$ be the number of subjects randomized to treatment 1 after $j$
assignments. Let $\boldX=(X_1,\ldots,X_n)$ be the responses based
on some primary outcome variable, and let $\mathbf{x}$ be the
realization. A valid test of the treatment effect can be conducted
permuting $\boldT$ in all possible ways [e.g., \citet{Leh86}, Chapter
5]. However, if one wishes to incorporate the
randomized design into the analysis, under restricted randomization, such
permutations are not equiprobable [e.g., \citet{RosLac02},
Chapter 7]. The family of \textit{linear rank tests} provides a large
class of test statistics with which to conduct randomization tests.
The form of the statistic is
$V(\mathbf{T})=\mathbf{a}_n^{\prime}\mathbf{T}$, for a
score vector $\mathbf{a}_n=(a_{1n}
-\bar{a}_n,\ldots,a_{nn}-\bar{a}_n)^{\prime}$, where $a_{jn}$ is some
function of the rank of the $j$th observation and
$\bar{a}_n=\sum_{j=1}^n a_{jn}/n$.
The $p$-value of the randomization test is computed with respect to
a reference set of sequences. The \textit{unconditional reference set}
contains all possible allocation sequences, including those that
give little or no information about the treatment effect (e.g.,
$1,1, \ldots, 1$). Also, the random numbers on each treatment arm,
$N_1(n)$ and~$N_2(n)$,
are ancillary statistics, and therefore the \textit{conditional reference
set} is preferred, which finds probabilities conditional on
$N_1(n)=n_1$, that is, the observed number on treatment $1$ [e.g., Cox
(\citeyear{Cox82}), \citet{Ber00}]. This leads to a~\textit
{conditional test}.

The literature
is largely silent on the subject of sequential monitoring of
randomization tests
(brief exceptions are found in Rosenberger and Lachin
[(\citeyear{RosLac02}), Section 7.10]
and \citet{ZhaRos08}, whose techniques only extend to one
interim inspection).
The computation of \textit{conditional} randomization tests
is also inherently difficult, even without sequential monitoring. We
address these issues in this paper by
proposing a technique, based on deriving exact conditional distributions
of randomization procedures, that leads to a simple computational
method for approximating the distribution of sequentially\vadjust{\goodbreak} computed
randomization tests. We also discuss the appropriate analog for
``information fraction'' in the context of a randomization model.
Our focus will be on one particular restricted randomization procedure,
Efron's (\citeyear{Efr71}) biased coin design, which induces a beautiful
closed-form combinatoric structure
that facilitates such an analysis. However, the technique can be applied
to any randomization procedure for which we can determine certain
exact conditional distributional results.

Let $\phi_{j+1}$ be a restricted randomization procedure such that
\[
\phi_{j+1}=\Pr\bigl(T_{j+1}=1|N_1(j)\bigr).
\]
Efron's (\citeyear{Efr71})
biased coin design is a restricted randomization procedure for
clinical trials that has exceptional properties: it balances
treatment assignments throughout the course of the trial with low
variability [e.g., \citet{Ant08}], and it mitigates
selection and accidental biases [\citet{RosLac02}]. Then
the biased coin design $\operatorname{BCD}(p)$ for a parameter $p \in[1/2,1]$,
$q=1-p$, is defined as
%
%e1.1 #&#
%
\begin{equation}\label{phij}
\phi_{j+1} = \cases{
1/2, &\quad when $N_1(j) = j/2$,\cr
p, &\quad when $N_1(j) < j/2, j = 0,1,2,\ldots,$ \cr
q, &\quad when $N_1(j) > j/2$.}
\end{equation}
Note that $p=0.5$ results in complete randomization and $p=1$
results in permuted blocks with block size 2. When $p<1$, the
design is \textit{fully randomized}, in that each subject will be
assigned to treatment randomly, which differs markedly from the
permuted block design, where some subjects in the tail of each block
are assigned with probability 1. Let $D_j=2N_1(j)-j$ be the
difference in numbers assigned to treatments 1 and 2;
$\{|D_{n}|\}_{n=1}^{\infty}$ forms an asymmetric random walk when
$p\in(1/2,1)$. \citet{MarRos10} derive the exact
distribution of $D_j$ for the $\operatorname{BCD}(p)$, from which the exact
distribution of $N_1(n)$ follows immediately:
%
%e1.2 #&#
%
\begin{eqnarray}\label{Markaryan}\qquad
&&P\bigl(N_1(n)=n_1\bigr)\nonumber\\[-8pt]\\[-8pt]
&&\qquad=
\cases{
\displaystyle p^{n_1} \sum_{l=0}^{n_1-1}\frac{n-2l}{n+2l}
\pmatrix{n_1+l\cr l}q^l,&\quad$\displaystyle n_1=\frac{n}{2}$,\vspace*{2pt}\cr
\displaystyle \frac{p^{n_1}}{2} \sum_{l=0}^{n_1}\frac
{n-n_1-l}{n-n_1+l} \pmatrix{n-n_1+l\cr l}q^{n-2n_1+l-1},&\quad
$\displaystyle 0\leq n_1 <
\frac{n}{2}$,\vspace*{2pt}\cr
\displaystyle \frac{p^{n-n_1}}{2} \sum_{l=0}^{n-n_1}\frac
{n_1-l}{n_1+l} \pmatrix{n_1+l\cr l}q^{2n_1-n+l-1},&\quad$\displaystyle \frac{n}{2}<
n_1 \leq n$.}\nonumber
\end{eqnarray}
Their paper also provides the exact expression for the
variance--covariance matrix of the treatment assignments
$\mathbf{T}$.

In this paper we provide the exact conditional distribution of
$N_1(n)$ given $N_1(j), 1\leq j < n$, and the expression for the
variance--covariance matrix of $\mathbf{T}$ given $N_1(n)$,
$\bolds{\Sigma}_{|n_1}$. While these are heretofore unknown results
on theoretical properties of a random walk, our primary interest is
that these results give us a computational method to approximate
conditional randomization tests following the $\operatorname{BCD}(p)$.
We then extend these results to the case where sequential analysis is
implemented in the course of a clinical trial.

\citet{RosLac02} distinguish among three techniques
that can be used to compute randomization tests: exact, Monte
Carlo and asymptotic. Exact tests are computationally intensive,
even with today's computing, and require networking algorithms
[\citet{autokey9}]. \citet{HolPen88}
developed a clever recursive algorithm to
determine the exact distribution of both conditional and
unconditional randomization tests following Efron's biased coin
design and applied it to a sample of size of $n=37$. It can be
assumed that such computational techniques would be able to solve
much larger problems with today's computing resources. While authors
have determined the asymptotic normality of conditional
randomization tests under various score functions and randomization
procedures [e.g., \citet{Smy88}], Efron's biased coin induces a
stationary distribution, and hence the test statistic may not be
asymptotically normal. This phenomenon was noted in a number of
papers, first by counterexample in \citet{SmyWei83} for the
unconditional test, and then by simulation by
\citet{HolPen88} for the conditional test.

\citet{MehPatSen88} use importance sampling to
estimate the conditional randomization test's $p$-value; their
technique employs an elegant, but complex, networking algorithm. The
efficiency of the estimator relies on the convergence to normality of
the test statistic, which may not hold under the biased coin design.
One might be able to modify the network algorithm in
\citet{autokey9} or the recursive algorithm in
\citet{HolPen88} to compute the exact distribution of sequentially
monitored conditional randomization tests, but here we provide a method
that is not very computationally intensive and allows us to sample
directly from the conditional reference set under a broad class of
restricted randomization procedures.

The paper is organized as follows. In Section \ref{sec2}, we present a method
for sampling directly from the conditional
distribution of $V(\mathbf{T})$, which facilitates the computation
of conditional tests.
We need to compute certain exact conditional probabilities to apply
this method, and we do
this for Efron's biased coin design.
We extend this application to
develop a computational technique for sequential monitoring of
conditional randomization tests in Section \ref{sec3}. In Section \ref
{sec4}, we
describe the analog of ``information'' in the context of a randomization
model. In defining a randomization-based
information fraction, we must derive the conditional
variance--covariance matrix of $\mathbf{T}$, which we
do for Efron's biased coin design. We draw conclusions in Section
\ref{sec5}. All the major proofs, some of which require careful
combinatorics, are relegated to the online supplement.

%s2 #&#
\section{Computation of conditional randomization tests}\label{sec2}

%s2.1 #&#
\subsection{Generating sequences from the unconditional reference
set}\label{sec21}
Suppose a total of $n$ subjects are randomized to two treatments. Let
$n_1$ be the observed number of assignments
on treatment $1$. A conditional randomization test can be approximated
by sampling sufficiently many sequences
from the conditional reference set, $\Omega_c$, the collection of
sequences that satisfy $N_1(n)=n_1$. This can be achieved by
generating sequences from the unconditional reference set, the set
of all possible assignments, and retaining those that belong to
$\Omega_c$.

Suppose at least $N_c$ number of sequences that satisfy
$N_1(n)=n_{1}$ are sufficient to approximate the conditional
randomization distribution of $V(\mathbf{T})$. Let $K$ sequences
be sampled, $\mathbf{T}_1,\ldots,\mathbf{T}_{K}$,
independently and with replacement from the unconditional reference
using $\phi_{j+1}$ as the sampling mechanism. This number of Monte
Carlo sequences must be large enough such that at least $N_c$
sequences satisfy the condition $N_1(n)=n_{1}$. The requisite number
of sequences, $K$, follows a negative binomial random variable with
parameters $\pi=P(N_{1}(n)=n_{1})$ and $r=N_c$ [\citet{ZhaRos10}].
Let $N$ denote a value in the range of $K$,
$N=N_c, N_c+1,\ldots.$ For $k=1,\ldots,N$, a sequence
$\mathbf{T}_k=\mathbf{t}$ is sampled from the unconditional
reference set with probability
%
%e2.1 #&#
%
\begin{equation}\label{Ch2unormalized}
f(\mathbf{t})=(1/2)
\prod_{j=1}^{n-1}(\phi_{j+1})^{t_{j+1}}(1-\phi_{j+1})^{1-t_{j+1}},
\end{equation}
where $t_{j+1}$ are the observed values of $T_{j+1}$. The $j$th
sampled sequence induces two Bernoulli random variables
\[
Y_j=
\cases{
1, &\quad if $N_{1}(n)=n_1$,\cr
0, &\quad otherwise,}
\]
and
\[
X_j=
\cases{
1, &\quad if $N_{1}(n)=n_1$ and $V(\mathbf{T}_j)\geq
v^*$,\cr
0, &\quad otherwise,}
\]
where $v^*$ is the observed value of the statistic. A strongly
consistent estimator for the $p$-value of the upper-tailed
conditional test can be computed as
%
%e2.2 #&#
%
\begin{equation}\label{Ch2776}
\hat{p}_c=\frac{\sum_{j=1}^{N}X_j}{\sum_{j=1}^{N}Y_j}.
\end{equation}

Table \ref{table1} reports the $95$th percentile of $K$ when sampling from
the unconditional reference set and $N_c=2500$ for Efron's biased coin
design. These sample sizes
%
%t1 #&#
%
\begin{table}
\caption{Approximate $95$th percentile of $K$ for various $n$,
$n_1$, $N_c=2500$}\label{table1}
\begin{tabular*}{\tablewidth}{@{\extracolsep{\fill}}lccc@{}}
\hline
$\bolds{n}$& $\bolds{n_1=0.45n}$ & $\bolds{n_1=0.48n}$ & $\bolds{n_1=0.50n}$\\
% heading
\hline
\multicolumn{4}{@{}c@{}}{$\operatorname{BCD}(p=2/3$)} \\
[4pt]
100 & 3,531,344 & 55,060 & 5117 \\
200 & 3,611,280,266 & 881,557 & 5117 \\
500 & $3\mbox{,}877\mbox{,}310\times10^{12}$ & 3,611,026,232 & 5117 \\
[4pt]
\multicolumn{4}{@{}c@{}}{$\operatorname{BCD}(p=3/4)$} \\
[4pt]
100 &114,384,212 & 156,865 & 3822 \\
200 & $6\mbox{,}754\mbox{,}269 \times10^6$ & 12,709,307 & 3822 \\
500 & $1\mbox{,}390\mbox{,}644 \times10^{21}$ & $6\mbox{,}754\mbox{,}269\times10^6$ & 3822 \\
\hline
\end{tabular*}
\vspace*{12pt}
\end{table}
are reasonable when there is perfect balance in the assignments,\vadjust{\goodbreak} but
increase considerably in the presence of imbalance. This technique
cannot be used in the presence of any imbalance.

%s2.2 #&#
\subsection{Our method: Generating sequences from the conditional
reference set}\label{sec22}

Rather than sampling too many sequences and discarding those that do
not satisfy the condition $N_1(n)=n_{1}$, it is more efficient to
sample directly from~$\Omega_c$---the collection of all
randomization sequences that satisfy the condition $N_1(n)=n_{1}$.
The set $\Omega_c$ will be called the conditional reference set. Let
$N_c$ randomization sequences, $\mathbf{T}_{1},
\ldots,\mathbf{T}_{N_c}$, be sampled independently and with
replacement strictly from $\Omega_c$, each with respective
probabilities $h(\mathbf{t}_{1}),\ldots,h(\mathbf{t}_{N_c})$.
For an upper-tailed test, the $k$th sampled sequence induces a
Bernoulli random variable
%
%e2.3 #&#
%
\begin{equation}\label{Ch21}
V_k=
\cases{
1, &\quad if $V(\mathbf{T}_k) \geq v^*$,\cr
0, &\quad otherwise.}
\end{equation}
The Monte Carlo estimator of the upper-tailed test's $p$-value is the
strongly consistent and unbiased estimator
$\bar{V}=\sum_{k=1}^{N_c}V_k/N_c$. (It may be possible to find
an estimator with a smaller variance, but we do not address
the issue of estimation of $p$-values in this paper.)

To guarantee a sequence from $\Omega_c$, $T_{j+1}$ in $\phi_{j+1}$
must be conditioned on both $N_1(j)$ and $N_1(n)$. Consequently, for
$0 \le m_j \le j$, the procedure
%
%e2.4 #&#
%
\begin{equation}\label{Ch23}
p_{j+1}=
\cases{
P\bigl(T_{j+1}=1|N_1(j)=m_j, N_1(n)=n_{1}\bigr), &\quad$1 \leq j \leq n-1$,\vspace*{2pt}\cr
P\bigl(T_{j+1}=1|N_1(n)=n_{1}\bigr), &\quad$j=0$,}\hspace*{-30pt}
\end{equation}
must be applied to generate a random sequence strictly from
$\Omega_c$. We now provide a general formula relating the
conditional and the unconditional reference sets, which facilitates the
generation of sequences from the conditional reference set for any
restricted randomization procedure of the form
$\phi_{j+1}=\Pr(T_{j+1}=1|N_1(j))$.
%
%th2.1 #&#
%
\begin{theorem}
For $n=1,2,3,\ldots,0 \leq n_1 \leq n$, $0 \leq j <
n$, $0\le m_j
\le j$ and $\phi_{j+1}(m_j)=P(T_{j+1}=1|N_1(j)=m_j)$, the rule
%
%e2.5 #&#
%
\begin{equation}\label{Ch24}
p_{j+1}=\cases{
\displaystyle \phi_{j+1}(m_j)\frac{P(N_1(n)=n_{1}| N_1(j+1)=m_j+1)}{P
(N_1(n)=n_{1}|N_1(j)=m_j)}, \vspace*{2pt}\cr
\hspace*{162pt}\qquad1 \leq j \leq
n-1,\vspace*{2pt}\cr
\displaystyle \phi_{j+1}(m_j)\frac{P(N_1(n)=n_{1}|T_{j+1}=1)}{P
(N_1(n)=n_{1})},
\qquad j=0,}\hspace*{-25pt}
\end{equation}
can be used to sample a sequence that satisfies $N_{1}(n)=n_{1}$.
\end{theorem}
\begin{pf}
The result follows from an application of Bayes theorem to
(\ref{Ch23}) and the Markovian property of $N_1$.
\end{pf}

Furthermore, for $k=1,\ldots,N_c$, a sequence
$\mathbf{T}_k=\mathbf{t}$ is sampled from $\Omega_c$ with
probability
%
%e2.6 #&#
%
\begin{equation}\label{Ch25}
h(\mathbf{t})=\prod_{j=0}^{n-1}(
p_{j+1})^{t_{j+1}}(1-p_{j+1})^{1-t_{j+1}}.
\end{equation}
In the simplest case, complete randomization,
$p_{j+1}=(n_1-m_j)/(n-j), 0 \leq j \leq n-1$, and this is the
random allocation rule [see \citet{RosLac02}], which is
sometimes used to fill permuted blocks.

The following theorem gives these probabilities for Efron's biased coin
design. The distribution of $N_1(n)$ given $N_1(j)=m_j$, $0 \leq m_j
\leq
j$, has three cases depending on the value of $m_j$ with respect to
$j$, $1\leq j < n$. Within each case, $P(N_1(n)=n_{1}|N_1(j)=m_j)$
depends the value of $n_1$ with respect to $n$, $j$ and $m_j$.
%
%th2.2 #&#
%
\begin{theorem}\label{th2.2}
Let $n=2,3,4,\ldots,1 \leq j < n$, $0 \leq m_j \leq j$ and $m_j \leq
n_1 \leq n-j+m_j$. Denote
\[
C(x,l) := \frac{x-l}{x+l}\pmatrix{x+l\cr l}
\quad\mbox{and}\quad D := \pmatrix{n-j\cr n_1-m_j}-\pmatrix{n-j\cr
n_1-j+m_j}.
\]
For the $\operatorname{BCD}(p)$:
\begin{longlist}[(1)]
\item[(1)] When $0\leq
m_j<j/2$, $P(N_1(n)=n_{1}|N_1(j)=m_j)$ is
\begin{eqnarray}
&\displaystyle \pmatrix{n-j\cr n_1-m_j}p^{n_1-m_j}q^{n-j-n_1+m_j}
\qquad\mbox{if } m_j \leq n_1< j-m_j, &\nonumber\\
&\displaystyle0.5p^{n_1-m_j}\sum_{l=0}^{n_1+m_j-j}C(n-n_1-m_j,l)q^{n-2n_1-1+l}& \nonumber\\
&\qquad\hspace*{8pt}{}+Dp^{n_1-m_j}q^{n-j-n_1+m_j}\qquad \mbox{if } j-m_j \leq
n_1 < n/2, &\nonumber\\
&\displaystyle p^{n_1-m_j}\sum_{l=0}^{n-j-n_1+m_j}C(n_1-m_j,l)q^l
\qquad\mbox{if } n_1=n/2,&\nonumber\\
&\displaystyle 0.5p^{n-n_1-m_j}\sum
_{l=0}^{n-j-n_1+m_j}C(n_1-m_j,l)q^{2n_1-n-1+l}&\nonumber\\
\eqntext{\mbox{if } n/2< n_1
\leq n-j+m_j.}
\end{eqnarray}

\item[(2)] When $m_j=j/2$,
\begin{eqnarray*}
&&P\bigl(N_1(n)=n_{1}|N_1(j)=m_j\bigr)\\
&&\qquad=P\bigl(N_1(n-j)=n_{1}-m_j\bigr),\qquad m_j \leq
n_1\leq n-j+m_j,
\end{eqnarray*}
where the unconditional distribution is derived in Markaryan and Rosenberger
(\citeyear{MarRos10}).

\item[(3)] When $j/2<m_j\leq j$,
$P(N_1(n)=n_{1}|N_1(j)=m_j)$ is
\begin{eqnarray*}
&\displaystyle 0.5p^{n_1+m_j-j}\sum
_{l=0}^{n_1-m_j}C(n-j-n_1+m_j,l)q^{n-2n_1-1+l}\qquad \mbox{if } m_j \leq
n_1 < n/2,& \\
&\displaystyle p^{n-j-n_1+m_j}\sum_{l=0}^{n_1-m_j}C(n-j-n_1+m_j,l)q^l
\qquad\mbox{if } n_1 = n/2,&\\
&\displaystyle0.5p^{n-j-n_1+m_j}\sum
_{l=0}^{n-n_1-m_j}{C}(n_1+m_j-j,l)q^{2n_1-n-1+l}&  \\
&\hspace*{7pt}{}+Dp^{n-j-n_1+m_j}q^{n_1-m_j}\qquad \mbox{if }
n/2 < n_1 \leq n-m_j,&\\
&\pmatrix{n-j\cr n_1-m_j}p^{n-j-n_1+m_j}q^{n_1-m_j}
\qquad\mbox{if } n-m_j < n_1 \leq n-j+m_j.&
\end{eqnarray*}
\end{longlist}
\end{theorem}
\begin{pf}
See Appendix A in the supplementary material [\citet{PlaRos}].
\end{pf}

Note that if $n=j$ and $n_1=m_j$, $P(N_1(n)=n_1|N_1(j)=m_j)=1$, and
if $m_j>n_1$ or $n-j< n_1-m_j$, $P(N_1(n)=n_1|N_1(j)=m_j)=0$. Also,
$P(N_1(n)=n_1|N_1(0)=0)=P(N_1(n)=n_1)$.

The procedure then follows by simply
generating $N_c$ sequences using (\ref{Ch24}). This allows us to
reduce the magnitude of the problem from the astronomical numbers in
Table \ref{table1} to just $N_c$.
A satisfactory value for $N_c$ can be obtained from the constraint
$\operatorname{MSE}(\bar{V})=p_c(1-p_c)/N_c\le1/4N_c\le
\varepsilon
$. For\vspace*{1pt}
$\varepsilon=0.0001$, $N_c\ge2500$. Higher precision in estimation is
possible by finding $N_c$ that ensures $P(|\bar{V}-p_c
|\leq0.1p_c)=0.99$, for instance. It follows that
$N_c\approx(2.576/0.1)^2(1-p_c)/p_c$. Thus, to estimate a $p$-value
as large as 0.04 with an error of $10\%$ of 0.04 with 0.99
probability, the Monte Carlo sample size must be $N_c=15\mbox{,}924$. If a
smaller $p$-value is expected, $N_c$ will be larger.

%t2 #&#
%
\begin{table}
\caption{Approximations for the upper $0.1$ tail of the
randomization distribution of the linear rank statistic by sampling
from $\Omega_c$; $N_c=2500$, $\operatorname{BCD}(0.6)$}\label{table2}
\begin{tabular*}{\tablewidth}{@{\extracolsep{\fill}}l c c c c@{}}
\hline
& $\bolds{n}$ & $\bolds{n_1}$ & \textbf{Exact} & \textbf{1000 Monte Carlo runs; mean (SD)}\\
\hline
$P(V(\boldT)\geq21.5)$ & \hphantom{0}30 & \hphantom{0}15 & 0.1057 & 0.1053 (0.0061) \\
$P(V(\boldT)\geq23)$ & \hphantom{0}30 & \hphantom{0}12 & 0.1009 & 0.1008 (0.0059) \\
$P(V(\boldT)\geq31)$ & \hphantom{0}40 & \hphantom{0}20 & 0.1011 & 0.1009 (0.0061) \\
$P(V(\boldT)\geq34)$ & \hphantom{0}40 & \hphantom{0}16 & 0.1000 & 0.0997 (0.0060) \\
$P(V(\boldT)\geq82)$ & 100 & \hphantom{0}50 & & 0.1055 (0.0060) \\
$P(V(\boldT)\geq113)$ & 100 & \hphantom{0}40 & & 0.1043 (0.0062) \\
$P(V(\boldT)\geq299)$ & 500 & 250 & & 0.1104 (0.0063) \\
$P(V(\boldT)\geq1000)$ & 500 & 200 & &
0.1030 (0.0058) \\
\hline
\end{tabular*}
\end{table}

Table \ref{table2} provides approximations for the upper $0.1$ tail of the
linear rank statistic with simple rank scores under the BCD$(0.6)$
randomization. For small samples sizes, we also provide the exact
$p$-value for comparison purposes; Monte Carlo estimates are very
close with small variability. As expected, the variability of the
estimates does not change across different sample sizes $n$. The
computational complexity of the sampling scheme for the BCD is
invariant to the value of $p$. Comparing Tables \ref{table1} and \ref
{table2}, the
conditional distribution method reduces the Monte Carlo sample size
to a few thousand.

Following stratification on known covariates, the computation of a
stratified linear rank test based on the conditional randomization
distribution is straightforward by summing the stratum-specific
linear-rank test statistics over $I$ independent strata. Using the
methodology described in this section, a sequence is sampled
independently from the conditional reference of each stratum; the
linear-rank statistic is evaluated in each stratum and the
stratum-specific test statistics are summed. The process is repeated
$N_c$ times, and the stratified test's $p$-value is estimated by the
proportion of summed statistics as or more extreme than the one
observed.

%s3 #&#
\section{Extension to sequential monitoring}\label{sec3}

Suppose there are $L-1$ interim inspections of the data after $1
\leq r_1 < r_2<\cdots<r_{L-1}< r_L=n$ patients responded. Let
$0<t_1<t_2<\cdots<t_{L-1}<t_L=1$ be the corresponding information
fraction at those inspections. For conditional tests, let
$N_{1}(r_1),N_{1}(r_2),\ldots,\allowbreak N_{1}(r_{L-1})$, and $N_{1}(r_L)=N_1(n)$
be the sample sizes randomized to treatment~$1$ after inspections
$1,\ldots,L$ and let $n_{11},\ldots,n_{1(L-1)}$, and $n_{1L}=n_1$ be
realizations of these sample sizes. Let the linear-rank
randomization test statistic computed at each of the inspections be
given by $V_{r_l}=\sum_{j=1}^{r_l}(a_{jr_l} -
\bar{a}_{r_l})T_j=\mathbf{a}^{\prime}_{r_l}\mathbf
{T}^{(r_l)}, l=1,\ldots,L$.
Using the alpha-spending function approach [\citet{LanDeM83}],
let $\alpha^*(t),t\in[0,1]$, be a nondecreasing function such that
$\alpha^*(0)=0$ and $\alpha^*(1)=\alpha$, the significance level of
the one-sided test. One such function is
$\alpha^*(t)=2-2\Phi(z_{\alpha/2}/\sqrt{t}), 0 < t \leq1;$
$\alpha^*(0)=0$, where $\Phi$ is the standard normal distribution
function and $z_{\alpha/2}=\Phi^{-1}(1-\alpha/2)$ [\citet
{LanDeM83}, \citet{OBrFle79}]. Following \citet{ZhaRos08},
the upper-tailed, conditional randomization test with~$L$
interim looks involves finding $d_1,\ldots,d_L$ such that
%
%e3.1 #&#
%
\begin{equation}\label{seqcond}
\cases{
\displaystyle P\bigl(V_{r_1} > d_1|N_{1}(r_1)=n_{11}\bigr)=\alpha^*(t_1),\vspace*{2pt}\cr
\displaystyle P\bigl(V_{r_1} \le d_1, V_{r_2} >
d_2|N_{1}(r_1)=n_{11},N_{1}(r_2)=n_{12}\bigr)=\alpha^*(t_2)-\alpha
^*(t_1),\vspace*{2pt}\cr
\displaystyle P\Biggl(V_{r_1} \le d_1, V_{r_2} \le d_2,V_{r_3} > d_3\Big|\bigcap_{j=1}^3
N_{1}(r_j)=n_{1j}\Biggr)=\alpha^*(t_3)-\alpha^*(t_2),\cr
\vdots\vspace*{2pt}\cr
\displaystyle P\Biggl(V_{r_1} \le d_1,\ldots,V_L
> d_L\Big|\bigcap_{j=1}^{L} N_{1}(r_j)=n_{1j}\Biggr)=\alpha-
\alpha^*(t_{L-1}).}\hspace*{-35pt}
\end{equation}
The asymptotic joint normality of these conditional distributions
has not been shown, except in the case of $L=2$ under the
generalized biased coin design [\citet{ZhaRos08}].

We express (\ref{seqcond}) in terms of univariate conditional
distributions, which are much easier to sample from than the joint
distributions in (\ref{seqcond}).
%
%le3.1 #&#
%
\begin{lemma}\label{lem3.1}
The set of conditions (\ref{seqcond}) is equivalent to
%
%e3.2 #&#
%
\begin{equation}\label{seqcond1}
\cases{
P\bigl(V_{r_1} > d_1|N_{1}(r_1)=n_{11}\bigr)=\alpha^*(t_1),\vspace*{2pt}\cr
\displaystyle P\Biggl( V_{r_2} > d_2|V_{r_1} \le d_1, \bigcap_{j=1}^2\{
N_{1}(r_j)=n_{1j}\}\Biggr)=
\frac{\alpha^*(t_2)-\alpha^*(t_1)}{1-\alpha^*(t_1)},\vspace*{2pt}\cr
\displaystyle P\Biggl(V_{r_3} > d_3\Big| \bigcap_{j=1}^2 \{V_{r_j} \le d_j\},\bigcap
_{j=1}^3\{ N_{1}(r_j)=n_{1j}\} \Biggr)=
\frac{\alpha^*(t_3)-\alpha^*(t_2)}{1-\alpha^*(t_2)},\cr
\vdots\vspace*{2pt}\cr
\displaystyle P\Biggl(V_n > d_L \Big|\bigcap_{j=1}^{L-1} \{V_{r_j} \le
d_j\},\bigcap_{j=1}^{L}\{ N_{1}(r_j)=n_{1j}\}\Biggr) =
\frac{\alpha-
\alpha^*(t_{L-1})}{1-\alpha^*(t_{L-1})}.}\hspace*{-35pt}\vadjust{\goodbreak}
\end{equation}
\end{lemma}
\begin{pf}
See Appendix B in the supplementary material [\citet{PlaRos}].
\end{pf}

At each inspection $l$ in (\ref{seqcond1}), the conditional
reference set is the collection of all sequences satisfying
$\bigcap_{i=1}^l\{ N_{1}(r_i)=n_{1i}\}$. The following theorem can
be used to sample sequences from such sets.
%
%th3.1 #&#
%
\begin{theorem}\label{theo3.1}
Let $1\leq l \leq L$, $r_0, r_1, r_2,\ldots,r_l$ and $n_{10},
n_{11},\ldots,n_{1l}$ be defined as before, with $r_0=0$ and
$n_{10}=0$. Let $k=1,\ldots,l$. For $r_{k-1} \leq j < r_{k}$,
$n_{1(k-1)}\le m_j \le j$ and
$\phi_{j+1}(m_j)=P(T_{j+1}=1|N_1(j)=m_j)$, the rule
%
%e3.3 #&#
%
\begin{equation}\label{seqsamp}
\psi_{j+1}=\phi_{j+1}(m_j)\frac
{P(N_1(r_{k})=n_{1k}|N_1(j+1)=m_j+1)}{P(N_1(r_{k})=n_{1k}|N_1(j)=m_j)}
\end{equation}
can be used to sample a sequence that satisfies $\bigcap_{i=1}^l\{
N_{1}(r_i)=n_{1i}\}$.
\end{theorem}
\begin{pf}
See Appendix C in the supplementary material [\citet{PlaRos}].
\end{pf}

Note that equation (\ref{seqsamp}) reduces succinctly to the
expected $\psi_{j+1}= (n_{1k}-m_j)/(r_k-j)$ for complete
randomization, $l=1,\ldots,L$, $k=1,\ldots,l$, $r_{k-1} \leq j < r_{k}$
and $n_{1(k-1)}\le m_j \le j$. For the $\operatorname{BCD}(p)$ the numerator and
the denominator of $\psi_{j+1}$ must be evaluated according to
Theorem \ref{th2.2}. To obtain a~sequence from the reference set satisfying
$\bigcap_{i=1}^l\{ N_{1}(r_i)=n_{1i}\}$, the sampling must be done
in $k=1,\ldots,l$ steps as follows:
\begin{longlist}[(3)]
\item[(1)] At stage $k=1$, apply $\psi_{j+1}$ with $r_0 \le j < r_1$
to sample the first $r_1$ assignments.
\item[(2)] At stage $k=2$, apply $\psi_{j+1}$ with $r_1 \le j < r_2$
to sample the next $r_2-r_1$ assignments.
\item[(3)] At stage $3 \le k \le l$, apply $\psi_{j+1}$ with $r_{k-1}
\le j < r_k$ to sample the next $r_k-r_{k-1}$ assignments.
\end{longlist}

Suppose a sample of size $N_c$ (sequences) is sufficient to estimate
a distribution quantile using some quantile estimator. The Monte
Carlo algorithm that estimates the boundary $d_1,\ldots,d_L$ for an
$\alpha$-level, upper-tailed, conditional randomization test with
$L-1$ interim inspections is as follows:
\begin{longlist}[(3)]
\item[(1)] At stage 1, generate $N_c$ randomization sequences of $r_1$
assignments
from the reference set satisfying $N_{1}(r_1)=n_{11}$. Evaluate
$V_{r_1}$ for each sequence;
estimate $d_1$ using the nonparametric quantile estimator of
\citet
{CheLaz10} based on the values of $V_{r_1}$.

\item[(2)] At stage 2, generate $N_c/(1-\alpha^*(t_1))$ randomization
sequences of
$r_2$ assignments from the reference set satisfying $\bigcap_{i=1}^2\{
N_{1}(r_i)=n_{1i}\}$.
For each sequence, evaluate $V_{r_1}$
using the first $r_1$ of $r_2$ assignments only. Retain those
sequences that satisfy
$\{V_{r_1} \leq d_1\}$. Evaluate $V_{r_2}$ for each retained
sequence. Estimate $d_2$ using the quantile estimator of \citet
{CheLaz10} based on the values of $V_{r_2}$.

\item[(3)] At stage $3\leq l \leq L$, generate $N_c/\prod
_{i=1}^{l-1}(1-[\alpha^*(t_i)-\alpha^*(t_{i-1})]/[1-\alpha
^*(t_{i-1})])$ randomization sequences of
$r_l$ assignments from the reference set satisfying $\bigcap_{i=1}^l\{
N_{1}(r_i)=n_{1i}\}$. Note that $\alpha^*(t_0)=0$ and $\alpha
^*(t_L)=\alpha$. For each sequence, evaluate
$V_{r_1}, V_{r_2},\ldots,V_{r_{l-1}}$
using the first $r_1, r_2,\ldots,r_{l-1}$ assignments, respectively.
Retain those sequences that satisfy
$\bigcap_{i=1}^{l-1} \{V_{r_i} \leq d_i\}$. Evaluate $V_{r_l}$ for
each retained
sequence. Estimate $d_l$ using the quantile estimator of \citet
{CheLaz10} based on the values of $V_{r_l}$.
\end{longlist}
Requiring that
$N_c/\prod_{i=1}^{l-1}(1-[\alpha^*(t_i)-\alpha^*(t_{i-1})]/[1-\alpha
^*(t_{i-1})])$
randomization sequences be sampled at stage $l$ simply ensures that
at least $N_c$ sequences are used for the estimation of $d_l$ at
each stage $l$.

%s4 #&#
\section{Randomization-based information}\label{sec4}

Fisher's information is defined under a population model, and hence
it is not defined in the context of randomi\-zation-based inference.
However, since the
Fisher's information approximates the inverse of the asymptotic
variance of the test, it seems reasonable to define the
randomization-based analog of information as the ratio of the
variances [\citet{RosLac02}].
%
%e4.1 #&#
%
\begin{equation}\label{information}
t_l=\frac{\mathbf{a}^{\prime}_{r_l}\bolds{\Sigma
}_{|r_l}\mathbf{a}_{r_l}}{\mathbf{a}^{\prime}_n\bolds
{\Sigma}_{|n}\mathbf{a}_n},
\end{equation}
where
$\bolds{\Sigma}_{|r_l}=\operatorname{Var}(\mathbf
{T}^{(r_l)}|N_{1}(r_1)=n_{11},\ldots,N_{1}(r_l)=n_{1l})$.
This requires specification of $\bolds{\Sigma}_{|r_l}$ and
$\bolds{\Sigma}_{|n}$. We now derive these
for Efron's biased coin design. We begin with three lemmas:
%
%le4.1 #&#
%
\begin{lemma}\label{lem4.1}
Let $n=2,3,\ldots$ and $ 0 \leq n_1\leq n$. Let
$\phi_i(a)=P(T_i=1|N_1(i-1)=a)$ and
$f_{j-1,b}^{(i,a+1)}=P(N_1(j-1)=b|N_1(i)=a+1)$. For $1 \leq i < j
\leq n$,
\begin{eqnarray*}
&&E\bigl(T_iT_j|N_1(n)=n_1\bigr)\\
&&\qquad=\frac{
\sum_{a=0}^{i-1}\phi_i(a)P(N_1(i-1)=a)
\sum_{b=a+1}^{j-1}\phi_j(b)f_{j-1,b}^{(i,a+1)}
f_{n,n_1}^{(j,b+1)}}{P(N_1(n)=n_1)}.
\end{eqnarray*}
The conditional probabilities $f_{j-1,b}^{(i,a+1)}$ and
$f_{n,n_1}^{(j,b+1)}$ are given by Theorem \ref{th2.2}.
\end{lemma}
\begin{pf}
The result follows from an application of Bayes theorem to
$P(T_i=1,T_j=1|N_1(n)=n_1)$ and the Markovian property of
$N_1$.
\end{pf}

Given that we observe $N_1(n)=n_1$, we now derive the
variance--covariance matrix of $\boldT$, denoted by
$\bolds{\Sigma}_{|n_1}$.
%
%le4.2 #&#
%
\begin{lemma}\label{lem4.2}
Let $n=1,2,\ldots, 0 \leq n_1\leq n$,
$\vartheta_{i|n_1}=E(T_i|N_{1}(n)=n_{1})$ and
$\phi_i(a)=P(T_i=1|N_1(i-1)=a)$. For the $\operatorname{BCD}(p)$
\[
\vartheta_{i|n_1}=\frac{\sum_{a=0}^{i-1}P(N_1(i-1)=a)
\phi_i(a)P(N_1(n)=n_1|N_1(i)=a+1)}{P(N_1(n)=n_1)},\vadjust{\goodbreak}
\]
where
\[
\vartheta_{1|n_1}=1/2P\bigl(N_1(n)=n_1|N_1(1)=1\bigr)/P\bigl(N_1(n)=n_1\bigr).
\]
If
$i<j$, the $(i,j)$th entry of $\bolds{\Sigma}_{|n_1}$ is
\[
\sigma_{ij}=
\frac{ \sum_{a=0}^{i-1}\phi_i(a)P(N_1(i-1)=a)
\sum_{b=a+1}^{j-1}\phi_j(b)f_{j-1,b}^{(i,a+1)}
f_{n,n_1}^{(j,b+1)}}{P(N_1(n)=n_1)}- \vartheta_{i|n_1}\vartheta
_{j|n_1}.
\]
If $i=j$, the $(i,j)$th entry of $\bolds{\Sigma}_{|n_1}$ is
\[
\sigma_{ij}=\vartheta_{i|n_1}(1-\vartheta_{i|n_1}).
\]
\end{lemma}
\begin{pf}
The result follows from an application of Bayes theorem to
$P(T_i=1|N_1(n)=n_1)$, the Markovian property of $N_1$
and Lemma \ref{lem4.1}.
\end{pf}
%
%le4.3 #&#
%
\begin{lemma}\label{lem4.3}
Let $1\leq l \leq L$, $r_0, r_1, r_2,\ldots,r_l$ and $n_{10},
n_{11},\ldots,n_{1l}$ be defined as before, with $r_0=n_{10}=0$. Let
$\phi_i(a)=P(T_i=1|N_1(i-1)=a)$, $k=1,\ldots,l$, and
$f_{i-1,a}^{(r_{k-1},n_{1(k-1)}
)}=P(N_1(i-1)=a|N_1(r_{k-1})=n_{1(k-1)})$.
Denote $\vartheta_{i|r_l}=E(T_i|\bigcap_{q=1}^l\{
N_{1}(r_q)=n_{1q}\})$ and
$\lambda_{ij|r_l}=E(T_iT_j|\bigcap_{q=1}^l\{
N_{1}(r_q)=n_{1q}\})$.

For $1 \leq k \leq l$, $r_{k-1} < i \leq r_{k}$,
\[
\vartheta_{i|r_l}=\frac{\sum_{a=n_{1(k-1)}}^{i-1}\phi
_i(a)f_{i-1,a}^{(r_{k-1},n_{1(k-1)})}f_{r_k,n_{1k}}^{
(i,a+1)}
}{P(N_1(r_k)=n_{1k}|N_1(r_{k-1})=n_{1(k-1)})}.
\]

For $1 \leq k \leq l$ and $r_{k-1} < i<j \leq r_{k}$,
\[
\lambda_{ij|r_l} =\frac{
\sum_{a=n_{1(k-1)}}^{i-1}\phi_i(a)f_{i-1,a}^{
(r_{k-1},n_{1(k-1)})}
\sum_{b=n_{1(k-1)+1}}^{j-1}\phi
_j(b)f_{j-1,b}^{(i,a+1)}f_{r_k,n_{1k}}^{(j,b+1
)}}{f_{r_k,n_{1k}}^{(r_{k-1},n_{1(k-1)})}}.
\]

For all other $i,j$,
\[
\lambda_{ij|r_l}=
E\Biggl(T_i\Big|\bigcap_{q=1}^l\{
N_{1}(r_q)=n_{1q}\}\Biggr)E\Biggl(T_j\Big|\bigcap_{q=1}^l\{
N_{1}(r_q)=n_{1q}\}\Biggr).
\]
The probabilities
$f_{i-1,a}^{(r_{k-1},n_{1(k-1)})}$,
$f_{r_k,n_{1k}}^{(i,a+1)}$, $f_{j-1,b}^{(i,a+1)}$,
$f_{r_k,n_{1k}}^{(j,b+1)}$ and
$f_{r_k,n_{1k}}^{(r_{k-1},n_{1(k-1)})}$ are given by
Theorem \ref{th2.2}.
\end{lemma}
\begin{pf}
See Appendix D in the supplementary material [\citet{PlaRos}].
\end{pf}

Finally, the closed form of $\bolds{\Sigma}_{|r_l}$ is given in
the following theorem, which follows immediately from Lemma \ref{lem4.3}:
%
%th4.1 #&#
%
\begin{theorem}\label{th4.1} Let $1\leq l \leq L$, $k=1,\ldots,l$,
$r_0, r_1,
r_2,\ldots,r_l$ and $n_{10},\break
n_{11},\ldots, n_{1l}$ be defined as before, with
$r_0=n_{10}=0$.\vadjust{\goodbreak}

The $(i,j)$th entry of $\bolds{\Sigma}_{|r_l}$ under the
$\operatorname{BCD}(p)$ is
\[
\sigma_{ij}=
\cases{
\lambda_{ij|r_l} - \vartheta_{i|r_l}\vartheta_{j|r_l}, &\quad if $i<j
\mbox{ and } r_{k-1}<i<j \leq
r_{k}$,\vspace*{2pt}\cr
\vartheta_{i|r_l}(1-\vartheta_{i|r_l}), &\quad if $i=j$,\vspace*{2pt}\cr
0, &\quad otherwise,}
\]
where $\vartheta_{i|r_l}$ and $\lambda_{ij|r_l}$ are given by Lemma
\ref{lem4.3}.
\end{theorem}

Although one can compute $\bolds{\Sigma}_{|n}$ and
$\bolds{\Sigma}_{|r_l}$ exactly using Theorem \ref{th4.1},
$\mathbf{a}^{\prime}_{n}$ in~(\ref{information}) remains unknown
at each interim inspection, since a portion of the data is
unobserved. One would have to interpolate sequentially the remaining
unknown data points in order to have a value for
$\mathbf{a}^{\prime}_n$ and an approximation for
(\ref{information}). Interpolating the unknown observations by
sampling with replacement the known observations is one way to
obtain a value for $\mathbf{a}^{\prime}_n$. In our simulations
with data generated from two normal distributions, $L=3$, $n=350$,
$n_1=174$ and assignments following the BCD$(3/4)$, the approximate
information fraction at the first interim look with $r_1=250$ and
$n_{11}=126$ was 0.3791, compared to the true information of 0.3759.
At the second interim look with $r_2=300$ and $n_{12}=148$, the
approximate information fraction was 0.6380, compared to the true
information of 0.6382.

We also simulate the probability of type I error in an example. For
this purpose, we generate a sample of $n=350$ observations from
$N(1,0.9)$ and simulate treatment assignments from BCD($p=3/4$). We
plan $L=3$ interim looks: at $r_1=250$, $r_2=300$ and $r_3=350$. The
observed number assigned to treatment 1 at each look was
$n_{11}=126$, $n_{12}=128$ and $n_{13}=174$. We compute the boundary
values using the algorithm in Section \ref{sec3}. Table \ref{table3}
gives the
estimated type I error rate ($\hat{\alpha}$) and standard
deviation over 1000 replications for this sequential conditional
test. The probability of type I error is preserved with low variability.

%t3 #&#
%
\begin{table}
\caption{Mean (SD) of
simulated $\alpha$ for an $\alpha= 0.05$ upper tail sequential test
over a Monte Carlo sample size of 1000, $N_c = 2500$, interpolating
the unknown observations by sampling with replacement}\label{table3}
\begin{tabular*}{\tablewidth}{@{\extracolsep{\fill}}l c c c c c c@{}}
\hline
\textbf{Look} $\bolds{l}$ & $\bolds{r_l}$
& $\bolds{n_{1l}}$ & $\bolds{t_l}$ & $\bolds{\alpha_l}$\tabnoteref{ta}
& $\bolds{\hat{d}_l}$ & $\bolds{\hat{\alpha}}$ \\
\hline
Look 1 & 250 & 126 & 0.3617 &0.0011 & 1709 & \\
Look 2 & 300 & 148 & 0.6248 &0.0121 & 1688 & \\
Look 3 & 350 & 174 & 1\hphantom{.0000} &0.0373 & 1501 & 0.0495 (0.0043) \\
\hline
\end{tabular*}
\tabnotetext[\mbox{$*$}]{ta}{$\alpha_l=\frac{\alpha^*(t_l)-\alpha^*(t_{l-1})}{1-\alpha^*(t_{l-1})}$.}
\end{table}

%s5 #&#
\section{Conclusions}\label{sec5}

We have provided a computational method to approximate conditional
randomization tests, which can be extended to clinical trials that
incorporate sequential\vadjust{\goodbreak} monitoring. The key is to determine certain
conditional probabilities from the particular randomization
procedure. These techniques apply to any restricted randomization
procedure of the form $\phi_{j+1}=\Pr(T_{j+1}=1|N_1(j))$ and for
which closed form conditional probabilities can be obtained. We have
derived the exact conditional distribution of $N_1(n)$, given
$N_1(j)$, for Efron's $\operatorname{BCD}(p)$ using combinatoric arguments, also
the conditional variance--covariance matrix of $\mathbf{T}$,
which allows computation of the information fraction.

The class of generalized biased coin designs (GBCD) [\citet
{Wei78}] does not have a known form for the exact conditional
distribution, and this remains an open problem. For the
sequential monitoring of conditional tests using the GBCD with one
interim look, \citet{ZhaRos08} derived the joint
asymptotic distribution of the interim and the final test
statistics, which allows for an asymptotic test.

\section*{Acknowledgments}

The authors thank Tigran Markaryan, Anindya Roy and the referees for
helpful comments.

\begin{supplement}[id=suppA]
\stitle{Supplement to ``Sequential monitoring with conditional
randomization tests''}
\slink[doi]{10.1214/11-AOS941SUPP} %[doi,text={...}] - jei reikia
%suskaldyti doi
\sdatatype{.pdf}
\sfilename{aos941\_supp.pdf}
\sdescription{The supplement contains Appendix~A (proof of Theorem
\ref{th2.2}),  Appendix B
(proof of Lemma \ref{lem3.1}),  Appendix C
(proof of Theorem \ref{theo3.1}), and Appendix
D (proof of Lemma \ref{lem4.3}).}
\end{supplement}

% imsref loaded by lrinkeviciute, 2012-01-12 09:02:34
% imsref loaded by lrinkeviciute, 2012-01-12 09:40:15
%

\printaddresses

\end{document}